\documentclass[leqno,a4paper,12pt]{article}
\usepackage{amsmath,amsthm}

\newtheorem{Thm}{\indent Theorem}[section]
\newtheorem{Prop}[Thm]{\indent Proposition}

\newtheorem{Cor}[Thm]{\indent Corollary}
\newtheorem{Var}[Thm]{\indent Variant}

\theoremstyle{definition}
\newtheorem{Def}[Thm]{\indent Definition}
\newtheorem{Rem}[Thm]{\indent Remark}

\newtheorem{Exs}[Thm]{\indent Examples}

\def\qed{{\hskip0pt\unskip\unskip\nobreak\hfil\penalty50
          \hskip1em\hbox{}\nobreak\hfil
          {\bf q.e.d.}%
          \parfillskip=0pt\finalhyphendemerits=0
          \par}\medskip}

\newenvironment{Proof}
               {{\it Proof.}\quad}
               {\qed}

\newenvironment{Proofof}[1]
               {{\it Proof of #1.}\quad}
               {\qed}

\newcommand{\Prime}{\kern3\fontdimen1\font$'$\kern-7\fontdimen1\font}

\long\def\forget#1{}

\long\def\beginSIDEREMARK#1\endSIDEREMARK
    {{\par\bigskip\advance\leftskip by 2cm
                  \advance\rightskip by -2cm\noindent
      {\bf Our own side remark:} #1
      \par\bigskip\noindent}}

\long\def\beginFORGET#1\endFORGET{#1}
\long\def\beginFORGET#1\endFORGET{}

\def\?{\ ???\ \immediate\write16{}%
\immediate\write16{Warning: There was still a question mark . . . }%
\immediate\write16{}}

\usepackage{amsmath}

\usepackage{exscale}
\usepackage{amssymb}

\newcommand{\BC}{{\mathbb{C}}}

\newcommand{\BQ}{{\mathbb{Q}}}

\newcommand{\BZ}{{\mathbb{Z}}}

\newcommand{\FH}{{\mathfrak{H}}}

\newcommand{\CA}{{\cal A}}
\newcommand{\CB}{{\cal B}}
\newcommand{\CC}{{\cal C}}

\newcommand{\CH}{{\cal H}}

\newcommand{\CK}{{\cal K}}

\newcommand{\Spec}{\mathop{{\bf Spec}}\nolimits}

\newcommand{\imm}{\mathop{{\rm im}}\nolimits}

\newcommand{\Ext}{\mathop{\rm Ext}\nolimits}

\newcommand{\Gr}{\mathop{\rm Gr}\nolimits}
\newcommand{\gr}{\mathop{\rm gr}\nolimits}
\newcommand{\Hom}{\mathop{\rm Hom}\nolimits}

\newcommand{\loccit}{[loc.$\;$cit.]}

\newbox\mybox
\def\arrover#1{\mathrel{
       \setbox\mybox=\hbox spread 1.4em{\hfil$\scriptstyle#1$\hfil}
       \vbox{\offinterlineskip\copy\mybox
             \hbox to\wd\mybox{\rightarrowfill}}}}
\def\larrover#1{\mathrel{
       \setbox\mybox=\hbox spread 1.4em{\hfil$\scriptstyle#1$\hfil}
       \vbox{\offinterlineskip\copy\mybox
             \hbox to\wd\mybox{\leftarrowfill}}}}

\def\ontoover#1{\mathrel{
       \setbox\mybox=\hbox spread 1.4em{\hfil$\scriptstyle#1$\hfil}
       \vbox{\offinterlineskip\copy\mybox
             \hbox to\wd\mybox{\rightarrowfill\hskip-2.8mm
                               $\rightarrow$}}}}
\def\leftontoover#1{\mathrel{
       \setbox\mybox=\hbox spread 1.4em{\hfil$\scriptstyle#1$\hfil}
       \vbox{\offinterlineskip\copy\mybox
             \hbox to\wd\mybox{$\leftarrow$\hskip-2.8mm
                               \leftarrowfill}}}}
\def\longto{\longrightarrow}
\def\into{\hookrightarrow}

\def\isoto{\arrover{\sim}}

\def\longinto{\lhook\joinrel\longrightarrow}

\usepackage[curve,matrix,arrow,cmtip]{xy}
\NoComputerModernTips

\def\myxymessage{\def\messagetext
   {Here an xy-pic diagram was omitted to speed up compilation . . . }
   \immediate\write16{\messagetext}
   \hbox{\bf \messagetext}}
\def\filxymatrix#1{\myxymessage}
\def\filxyarray#1{\myxymessage}

\newdir^{ (}{{}*!/-3pt/\dir^{(}}
\newdir_{ (}{{}*!/-3pt/\dir_{(}}
\newdir^{ )}{{}*!/+3pt/\dir^{)}}
\newdir_{ )}{{}*!/+3pt/\dir_{)}}

\def\rscript#1{\hbox to 0pt{$\scriptstyle#1$\hss}}

\let\oldbullet\bullet
\def\bullet{{\mathchoice{\oldbullet}%
                        {\oldbullet}%
                        {\scriptscriptstyle\oldbullet}%
                        {\oldbullet}}}

\newcommand{\CHeffM}{\mathop{CHM^{eff}(k)}\nolimits}
\newcommand{\CHM}{\mathop{CHM(k)}\nolimits}
\newcommand{\CHeffQM}{\mathop{\CHeffM_F}\nolimits}
\newcommand{\CHQM}{\mathop{\CHM_F}\nolimits}
\newcommand{\DeffgM}{\mathop{DM^{eff}_{gm}(k)}\nolimits}
\newcommand{\DeffQgM}{\mathop{\DeffgM_F}\nolimits}
\newcommand{\DgM}{\mathop{DM_{gm}(k)}\nolimits}
\newcommand{\DQgM}{\mathop{\DgM_F}\nolimits}
\newcommand{\QAM}{\mathop{M \! A(k)_F}\nolimits}
\newcommand{\DAM}{\mathop{DM \! A(k)}\nolimits}
\newcommand{\DQAM}{\mathop{\DAM_F}\nolimits}
\newcommand{\QDM}{\mathop{M \! D(k)_F}\nolimits}
\newcommand{\QATM}{\mathop{M \! AT(k)_F}\nolimits}
\newcommand{\DQATM}{\mathop{D \! \QATM}\nolimits}

\newcommand{\QDTM}{\mathop{M \! DT(k)_F}\nolimits}
\newcommand{\DQDTM}{\mathop{D \! \QDTM}\nolimits}
\newcommand{\DTM}{\mathop{DMT(k)}\nolimits}
\newcommand{\DQTM}{\mathop{\DTM_F}\nolimits}

\newcommand{\ShN}{\mathop{Shv_{Nis}(SmCor(k))}\nolimits}
\newcommand{\Mgm}{\mathop{M_{gm}}\nolimits}

\begin{document}

%

\hfuzz=3pt
\overfullrule=10pt                   


\setlength{\abovedisplayskip}{6.0pt plus 3.0pt}
\setlength{\belowdisplayskip}{6.0pt plus 3.0pt}
\setlength{\abovedisplayshortskip}{6.0pt plus 3.0pt}
\setlength{\belowdisplayshortskip}{6.0pt plus 3.0pt}

\setlength{\baselineskip}{13.0pt}
\setlength{\lineskip}{0.0pt}
\setlength{\lineskiplimit}{0.0pt}

%
%

\title{Notes on Artin--Tate motives
\forget{
\footnotemark
\footnotetext{To appear in ....}
}
}
\author{\footnotesize by\\ \\
\mbox{\hskip-2cm
\begin{minipage}{6cm} \begin{center} \begin{tabular}{c}
J\"org Wildeshaus \footnote{
Partially supported by the \emph{Agence Nationale de la
Recherche}, project no.\ ANR-07-BLAN-0142 ``M\'ethodes \`a la
Voevodsky, motifs mixtes et G\'eom\'etrie d'Arakelov''. }\\[0.2cm]
\footnotesize LAGA\\[-3pt]
\footnotesize UMR~7539\\[-3pt]
\footnotesize Institut Galil\'ee\\[-3pt]
\footnotesize Universit\'e Paris 13\\[-3pt]
\footnotesize Avenue Jean-Baptiste Cl\'ement\\[-3pt]
\footnotesize F-93430 Villetaneuse\\[-3pt]
\footnotesize France\\
{\footnotesize \tt wildesh@math.univ-paris13.fr}
\end{tabular} \end{center} \end{minipage}
\hskip-2cm}
}
\date{January 14, 2010}
\maketitle
\begin{abstract}
\noindent In this paper,
we study the main structural properties
of the triangulated category of Artin--Tate motives
over a perfect base field $k$. We first analyze its
weight structure, 
building on the main results of \cite{Bo}.
We then study its $t$-structure,
when $k$ is algebraic over $\BQ$,
gene\-ralizing the main result of \cite{L}.
We finally exhibit the interaction of the weight
and the $t$-structure. When $k$ is a number field,
this will give a useful criterion
identifying the weight structure \emph{via} realizations. \\

\noindent Keywords: Artin--Tate motives, Dirichlet--Tate motives,
weight structures, $t$-structures, realizations.

\end{abstract}


\bigskip
\bigskip
\bigskip

\noindent {\footnotesize Math.\ Subj.\ Class.\ (2000) numbers: 
14F42 (14C35, 14D10, 19E15, 19F27). }

\eject

\tableofcontents

\bigskip


%
%

\setcounter{section}{-1}
\section{Introduction}
\label{Intro}



The aim of this article is to exhibit the basic structural properties
of the \emph{triangulated category of Artin--Tate motives}
over a fixed perfect base field $k$. The definition of this category
will be recalled, and a number
of generalizations will be defined in Section~\ref{1}.
Roughly speaking, the properties we shall be interested in, then
fall into two classes. \\
 
First (Section~\ref{2}), we apply the main results of \cite{Bo}
to Artin--Tate motives. More precisely,
we show (Theorem~\ref{2E}) that the \emph{weight structure}
of \loccit , defined on the category of \emph{geometrical motives}
\cite{V}, induces a weight structure on the 
triangulated category of Artin--Tate motives.
We also give a very explicit description of the \emph{heart}
of the latter,
showing in particular that it is Abelian semi-simple. \\

Second (Section~\ref{3}), we generalize the main result from \cite{L}
from Tate motives to Artin--Tate motives,
when the base field is algebraic over $\BQ$.
More precisely, we show (Theorem~\ref{1Main})
that under this hypothesis, there is
a non-degenerate \emph{$t$-structure} on
the triangulated category of Artin--Tate motives. 
The strategy of proof is identical to the one used by Levine. 
Using the main result of \cite{W4}, we show (Corollary~\ref{1N}) 
that this triangulated category is canonically equivalent
to the bounded derived category of its \emph{heart} (formed with respect
to the $t$-structure), \emph{i.e.}, of the Abelian category of \emph{mixed Artin--Tate
motives}. \\

Our main interest lies then in the simultaneous application of both
points of view: that of weight structures and that of $t$-structures. 
Still assuming that $k$ is algebraic over $\BQ$,
we give a characterization  
(Theorem~\ref{2Main})
of the weight structure on the triangulated category of Artin--Tate motives
in terms of the $t$-structure. 
Specializing further to the case of number fields,
we get a powerful criterion (Theorem~\ref{2H}), 
allowing to identify
the weight structure \emph{via} the Hodge theoretic or
$\ell$-adic \emph{realization}. \\

We should warn the reader that 
all our constructions are \emph{a priori} with $\BQ$-coefficients.
This seems to be necessary for at least two reasons. First,
the \emph{triangulated category of Artin motives}
is not known to admit a $t$-structure;
by contrast, such a structure becomes obvious after 
tensoring with $\BQ$
(see Section~\ref{1}). This $t$-structure is
at the very basis of our construction. Second, as pointed out in \cite{L}, 
the existence of the $t$-structure on
the \emph{triangulated category of Tate motives} necessitates (and is in fact
equivalent to)
the validity of the Beilinson--Soul\'e vanishing conjecture; but this vanishing
is only known (for algebraic base fields) after tensoring with $\BQ \, $. \\

Part of this work was done while I was enjoying a 
\emph{modulation de service pour les porteurs de projets de recherche},
granted by the \emph{Universit{\'e} Paris~13}. 
I wish to thank P.~J\o rgensen, B.~Kahn and M.~Levine for 
useful comments. \\

{\bf Notation and conventions}: Throughout the article, 
$k$ denotes a fixed perfect base field. 
The notation of this paper follows
that of \cite{V}. We refer to \loccit \ for
the de\-fi\-nition of the triangulated 
categories $\DeffgM$ and $\DgM$ of (effective) geometrical
motives over $k$. Let $F$ be a commutative $\BQ$-algebra.
The notation $\DeffQgM$ and $\DQgM$ stands 
for the $F$-linear analogues of these triangulated categories
defined in \cite[Sect.~16.2.4
and Sect.~17.1.3]{A}. 
Similarly, let us denote by $\CHeffM$ and $\CHM$ the categories 
opposite to the categories of (effective) Chow motives, and by
$\CHeffQM$ and $\CHQM$ the pseudo-Abelian
completion of the category $\CHeffM \otimes_\BZ F$ and 
$\CHM \otimes_\BZ F$, respectively. 
Using \cite[Cor.~2]{V2} (\cite[Cor.~4.2.6]{V} if $k$ 
admits resolution of singularities), 
we canonically identify 
$\CHeffQM$ and $\CHQM$ with
a full additive sub-category of $\DeffQgM$ and $\DQgM$, respectively.


\bigskip


%
%

\section{Definition and first properties}
\label{1}



Fix a commutative $\BQ$-algebra $F$,
which we suppose to be semi-simple and Noetherian,
in other words, a finite direct product of fields
of characteristic zero.
In this section, we recall the definition of the 
$F$-linear triangulated category 
of Artin--Tate motives (Definition~\ref{1C}), and define
a number of variants, indexed by certain 
sub-categories of the category of discrete representations 
of the absolute Galois group of our perfect base field $k$
(Definition~\ref{1F}). We then start the analysis of this category,
following the part of \cite{L} valid without additional assumptions
on $k$.  \\

For any integer $m$, there is defined a Tate object $\BZ(m)$ in $\DgM$,
which belongs to $\DeffgM$ if $m \ge 0$ \cite[p.~192]{V}. 
We shall use the same notation when we consider $\BZ(m)$ as an object
of $\DQgM$.

\begin{Def}[cmp.\ {\cite[Def.~3.1]{L}}] \label{1A}
Define the \emph{triangulated 
category of Tate motives over $k$}
as the strict full triangulated sub-category $\DQTM$ of $\DQgM$ generated by
the $\BZ(m)$, for $m \in \BZ$. 
\end{Def}

Recall that by definition, a strict sub-category is closed under
isomorphisms in the ambient category.
It is easy to see that $\DQTM$
is tensor triangulated. 

\begin{Def} \label{1B}
Define the \emph{triangulated category of Artin motives over $k$}
as the pseudo-Abelian completion of the strict
full triangulated sub-category $\DQAM$ of $\DeffQgM$ generated by
the motives $\Mgm(X)$ of smooth zero-dimensional schemes $X$ over $k$. 
\end{Def}

This category is again tensor triangulated. 

\begin{Def} \label{1C}
Define the \emph{triangulated
category of Artin--Tate motives over $k$}
as the strict full tensor triangulated sub-category $\DQATM$ of $\DQgM$ 
generated by $\DQAM$ and $\DQTM$. 
\end{Def}

The following observation \cite[Remark~2 on p.~217]{V} is vital.

\begin{Prop} \label{1D}
The triangulated category $\DQAM$ of Artin motives is canonically equivalent
to $D^b(\QAM)$, the bounded derived category of the Abelian category
$\QAM$ of discrete representations of the absolute
Galois group of $k$ in finitely generated $F$-modules.
\end{Prop}  
 
More precisely, if $X$ is smooth and zero-dimensional over $k$,
and $\bar{k}$ a fixed algebraic closure of $k$, then the absolute Galois
group of $k$, when identified with the group of automorphisms of $\bar{k}$
over $k$, acts canonically on the set of $\bar{k}$-valued points of $X$.
The object of $\QAM$ corresponding to $M(X)$ under the equivalence
of Proposition~\ref{1D}
is nothing but the formal $F$-linear envelope of this set, with the
induced action of the Galois group.
Note that the category $\QAM$ is semi-simple.

\begin{Cor} \label{1E}
There is a canonical non-degenerate $t$-structure on 
the category $\DQAM$. Its heart is equivalent to $\QAM$.
\end{Cor}

By contrast \cite[Remark~1 on p.~217]{V}, 
it is not clear how to construct a non-degenerate $t$-structure
on the triangulated category $\DAM$ of \emph{zero motives}
(whose $F$-linearization equals $\DQAM$). \\

For the rest of this section, let us identify the triangulated categories
$\DQAM$ and $D^b(\QAM)$ \emph{via} the equivalence of Proposition~\ref{1D}.
Let us also 
fix a strict
full Abelian semi-simple $F$-linear tensor sub-category $\CA$ of $\QAM$,
contai\-ning the category $triv$
of objects of $\QAM$ on which the Galois group acts trivially. 

\begin{Def} \label{1F}
Define $D\CA T$ as the strict full tensor triangulated 
sub-category of $\DQATM$ generated by $\CA$,
and by $\DQTM$.
\end{Def}

\begin{Exs} 
(a)~When $\CA$ equals $\QAM$, then $D\CA T = \DQATM$. \\[0.1cm]
(b)~When $\CA$ equals $triv$, 
then $D\CA T = \DQTM$.
\end{Exs}

Let us agree to set $\BZ(n/2):=0$ for odd integers $n$. 
For any object $M$ of $D\CA T$ and any integer $n$, let us write $M(n/2)$
for the tensor product of $M$ and $\BZ(n/2)$. \\

Following \cite{L},
let us first define $D\CA T_{[a,b]}$ as the full triangulated sub-category
of $D\CA T$ generated by the objects $N(m)$, for $N \in \CA$
and $a \le -2m \le b$, for integers $a \le b$ (we allow $a = - \infty$
and $b = \infty$). We denote $D\CA T_{[a,a]}$ by $D\CA T_a$.

\begin{Prop} \label{1G}
The category $D\CA T_a$ is zero for $a \in \BZ$ odd. For $a \in \BZ$ even,
the exact functor
\[
D\CA T_a \longto \DQAM \; , \; 
M \longmapsto M(a/2)
\]
induces an equivalence between $D\CA T_a$ 
and the bounded derived category of $\CA$
(which is equal to the $\BZ$-graded category 
$\Gr_\BZ \CA = \oplus_{m \in \BZ} \CA$
over $\CA$).
\end{Prop} 

\begin{Proof}
By construction, the functor is exact, and identifies $D\CA T_a$
with the full triangulated sub-category of $\DQAM$ of objects,
whose cohomology lies in $\CA$. Recall that we identified
the categories $\DQAM$ and $D^b(\QAM)$. It remains to see that
the obvious exact functor
\[
D^b(\CA) \longto D^b(\QAM)
\]
is fully faithful. But this an immediate consequence of the fact
that the Abelian categories $\CA$ and $\QAM$ are semi-simple.
\end{Proof}

In particular, there is a canonical $t$-structure 
$(D\CA T_a^{\le 0} , D\CA T_a^{\ge 0})$
on $D\CA T_a$: the category $D\CA T_a^{\le 0}$ is the full sub-category
of $D\CA T_a$ generated by objects $N(-a/2)[r]$, for $N \in \CA$ and $r \ge 0$,
and $D\CA T_a^{\ge 0}$ is the full sub-category
gene\-rated by objects $N(-a/2)[r]$, for $N \in \CA$ and $r \le 0$. 
If $a$ is even, then the category $\CA$ is equivalent to 
the heart $\CA T_a$ of this canonical $t$-structure 
\emph{via} the functor $N \mapsto N(-a/2)$. \\

Second, we construct auxiliary $t$-structures.

\begin{Prop}[cmp.\ {\cite[Lemma~1.2]{L}}] \label{1H}
Let $a \le n \le b$. Then the pair $(D\CA T_{[a,n]},D\CA T_{[n+1,b]})$
defines a $t$-structure on $D\CA T_{[a,b]}$.
\end{Prop}

\begin{Proof}
Imitate the proof of \cite[Lemma~1.2]{L}. The decisive ingredient is 
the following generalization of the vanishing from \cite[Def.~1.1~i)]{L}:
\[
\Hom_{D\CA T} \bigl( N_1(m_1)[r],N_2(m_2)[s] \bigr) = 0 \; , \;
\forall \, m_1 > m_2 \, , \, N_1, N_2 \in \CA \, , \, r,s \in \BZ \; .
\]
It holds because $\Hom_{D\CA T} = \Hom_{\DQgM}$,
and $\Hom_{DM_{gm}(\bullet)_F}$ satisfies descent
for finite extensions $L / k$ of the base field. 
Choosing an extension $L$ splitting both 
$N_1$ and $N_2$ therefore allows to deduce the desired vanishing
from that of
\[
\Hom_{DM_{gm}(L)_\BQ} \bigl( \BZ(m_1)[r],\BZ(m_2)[s] \bigr) \; .
\]
\end{Proof}

\begin{Rem}
(a)~The above proof uses 
the relation of $K$-theory of $L$
tensored with $\BQ \,$, with $\Hom_{DM_{gm}(L)_\BQ}$.
This relation is established by work of Bloch \cite{Bl1,Bl2}
(see \cite[Section~II.3.6.6]{L2}), and will be used again in 
the proofs of Theorem~\ref{1Main} and Variant~\ref{1Var}. \\[0.1cm]
(b)~Levine pointed out that
the $t$-structures from Proposition~\ref{1H}, for varying $n$,
can be used to show that the category $D \CA T$ is pseudo-Abelian. 
We shall give an alternative proof of this result in Section~\ref{2},
using Bondarko's theory of weight structures (Corollary~\ref{2Ea}). 
\end{Rem}

Note that since $D\CA T_{[a,n]}$ and $D\CA T_{[n+1,b]}$ 
are themselves triangulated, the $t$-structure from Proposition~\ref{1H} 
is necessarily degenerate. 
As in \cite[Sect.~1]{L}, denote the truncation functors
by
\[
W_{\le n} : D\CA T_{[a,b]} \longto D\CA T_{[a,n]}
\]
and
\[
W_{\ge n+1} : D\CA T_{[a,b]} \longto D\CA T_{[n+1,b]} \; ,
\]
and note that for fixed $n$, they are compatible with change of $a$ or $b$.
Write $\gr_n$ for the composition of $W_{\le n}$ and $W_{\ge n}$
(in either sense). The target of this functor is the category $D\CA T_n$. 
We are now ready to set up the data necessary for the $t$-structure
we shall actually be interested in.

\begin{Def}[cmp.\ {\cite[Def.~1.4]{L}}] \label{1I}
Fix $a \le b$ (we allow $a = - \infty$ and $b = \infty$). \\[0.1cm]
(a)~Define $D\CA T_{[a,b]}^{\le 0}$ as the full sub-category of $D\CA T_{[a,b]}$
of objects $M$ such that $\gr_n \! M \in D\CA T_n^{\le 0}$
for all integers $n$ such that $a \le n \le b$. \\[0.1cm]
(b)~Define $D\CA T_{[a,b]}^{\ge 0}$ as the full sub-category of $D\CA T_{[a,b]}$
of objects $M$ such that $\gr_n \! M \in D\CA T_n^{\ge 0}$
for all integers $n$ such that $a \le n \le b$.
\end{Def}

As we shall see (Theorem~\ref{1Main}, Variant~\ref{1Var}),
the pair
$(D\CA T_{[a,b]}^{\le 0} , D\CA T_{[a,b]}^{\ge 0})$ defines a $t$-structure
on $D\CA T_{[a,b]}$, provided that the base field
$k$ is algebraic over $\BQ \,$. In particular,
we then get a canonical $t$-structure on $D\CA T$. 
The vital point will be the validity of the Beilinson--Soul\'e
vanishing conjecture for all finite field extensions of $k$.


\bigskip


%
%

\section{The motivic weight structure}
\label{2}



The purpose of this section
is to first review Bondarko's definition of weight structures
on triangulated categories, and his result on the existence of such
a weight structure on the categories $\DgM$ and $\DQgM$ \cite{Bo}. 
We shall then show 
(Theorem~\ref{2E}~(a)) that
the latter induces a weight structure on any of the triangulated categories
constructed in Section~\ref{1}. We shall also get a very explicit description of the heart
of this weight structure (Theorem~\ref{2E}~(b), (c)). 

\begin{Def}[cmp.\ {\cite[Def.~1.1.1]{Bo}}] \label{2A}
Let $\CC$ be a triangulated category. A \emph{weight structure on $\CC$}
is a pair $w = (\CC_{w \le 0} , \CC_{w \ge 0})$ of full 
sub-categories of $\CC$, such that, putting
\[
\CC_{w \le n} := \CC_{w \le 0}[n] \quad , \quad
\CC_{w \ge n} := \CC_{w \ge 0}[n] \quad \forall \; n \in \BZ \; ,
\]
the following conditions are satisfied.
\begin{enumerate}
\item[(1)] The categories
$\CC_{w \le 0}$ and $\CC_{w \ge 0}$ are 
Karoubi-closed (\emph{i.e.}, closed under retracts formed in $\CC$).
\item[(2)] (Semi-invariance with respect to shifts.)
We have the inclusions
\[
\CC_{w \le 0} \subset \CC_{w \le 1} \quad , \quad
\CC_{w \ge 0} \supset \CC_{w \ge 1}
\]
of full sub-categories of $\CC$.
\item[(3)] (Orthogonality.)
For any pair of objects $M \in \CC_{w \le 0}$ and $N \in \CC_{w \ge 1}$,
we have
\[
\Hom_{\CC}(M,N) = 0 \; .
\]
\item[(4)] (Weight filtration.)
For any object $M \in \CC$, there exists an exact triangle
\[
A \longto M \longto B \longto A[1]
\]
in $\CC$, such that $A \in \CC_{w \le 0}$ and $B \in \CC_{w \ge 1}$.
\end{enumerate}
\end{Def}

It is easy to see that for any integer $n$ and any object $M \in \CC$,
there is an exact triangle
\[
A \longto M \longto B \longto A[1]
\]
in $\CC$, such that $A \in \CC_{w \le n}$ and $B \in \CC_{w \ge n+1}$.
By a slight generalization of the terminology introduced in  
condition~\ref{2A}~(4), we shall refer to any such exact triangle
as a weight filtration of $M$.

\begin{Rem} \label{2B}
Our convention concerning the sign of the weight is opposite 
to the one from \cite[Def.~1.1.1]{Bo}, \emph{i.e.}, we exchanged the
roles of $\CC_{w \le 0}$ and $\CC_{w \ge 0}$. 
\end{Rem}

\begin{Def}[{\cite[Def.~1.2.1]{Bo}}] \label{2C}
Let $w = (\CC_{w \le 0} , \CC_{w \ge 0})$ be a weight structure on $\CC$.
The \emph{heart of $w$} is the full additive sub-category $\CC_{w = 0}$
of $\CC$ whose objects lie 
both in $\CC_{w \le 0}$ and in $\CC_{w \ge 0}$.
\end{Def}

One of the main results of \cite{Bo} is the following.

\begin{Thm}[{\cite[Sect.~6]{Bo}}] \label{2D}
(a)~If $k$ is of characteristic zero,
then there is a canonical weight structure on the category $\DeffgM$.
It is uniquely characterized by the requirement that its heart equal 
$\CHeffM$. \\[0.1cm]
(b)~If $k$ is of characteristic zero,
then there is a canonical weight structure on the category $\DgM$,
extending the weight structure from (a).
It is uniquely characterized by the requirement that its heart equal 
$\CHM$. \\[0.1cm]
(c)~Let $F$ be a commutative $\BQ$-algebra.
Analogues of statements (a) and (b) hold for the $F$-linearized
categories $\DeffQgM$, $\CHeffQM$, $\DQgM$, and $\CHQM$,
and for a perfect base field $k$ of arbitrary characteristic.
\end{Thm}

Let us refer to any of these weight structures as \emph{motivic}.
For a concise review of the main ingredients of Bondarko's proof, see 
\cite[Sect.~1]{W3}. \\

Now fix a finite direct product $F$ of fields of characteristic zero,
and a full Abelian $F$-linear tensor 
sub-category $\CA$ of $\QAM$,
contai\-ning the category $triv$. 
Recall (Definition~\ref{1F}) that 
\[
D\CA T \subset \DQATM \subset \DQgM 
\]
denotes the strict full tensor triangulated 
sub-category generated by $\CA$, and by the triangulated category $\DQTM$
of Tate motives. 
Intersecting with $D\CA T$,
the motivic weight structure $(\DgM_{F,w \le 0} , \DgM_{F,w \ge 0})$
from Theorem~\ref{2D}~(c) yields a pair 
\[
w := w_{\CA} := \bigl(D\CA T_{w \le 0} , D\CA T_{w \ge 0}\bigr)
\]
of full sub-categories of $D\CA T$. 

\begin{Thm} \label{2E}
(a)~The pair $w$ is a weight structure on $D\CA T$. \\[0.1cm]
(b)~The heart $D\CA T_{w = 0}$ equals the intersection of $D\CA T$ and $\CHQM$.
It ge\-ne\-rates the triangulated category $D\CA T$. 
It is Abelian semi-simple.
Its objects are finite direct sums of objects of the form
$N(m)[2m]$, for $N \in \CA$ and $m \in \BZ$. \\[0.1cm]
(c)~The functor from the $\BZ$-graded category $\Gr_\BZ \CA$
over $\CA$ to $D\CA T_{w = 0}$
\[
\Gr_\BZ \CA = \bigoplus_{m \in \BZ} \CA \longto D\CA T_{w = 0} \; , \;
(N_m)_{m \in \BZ} \longmapsto \oplus_{m \in \BZ} N_m(m)[2m]
\]
is an equivalence of categories.
\end{Thm}

\begin{Proof}
Define $\CK$ as the full additive sub-category of $D\CA T$ 
of objects, which are finite direct sums of objects of the form
$N(m)[2m]$, for $N \in \CA$ and $m \in \BZ$.
Note that $\CK$ generates the triangulated category
$D\CA T$. All objects of $\CK$
are Chow motives. In particular, by orthogonality~\ref{2A}~(3)
for the motivic weight structure (see \cite[Cor.~4.2.6]{V}), 
$\CK$ is \emph{negative}, \emph{i.e.},
\[
\Hom_{D\CA T}(M_1,M_2[i]) = \Hom_{\DQgM}(M_1,M_2[i]) = 0
\]
for any two objects $M_1, M_2$ of $\CK$, and any integer $i > 0$.
Therefore,
\cite[Thm.~4.3.2~II~1]{Bo} can be applied to ensure the existence
of a weight structure $v$ on $D\CA T$, uniquely characterized by the
property of containing $\CK$ in its heart. Furthermore
\cite[Thm.~4.3.2~II~2]{Bo}, the heart $D\CA T_{v = 0}$ of $v$
is equal to the category $\CK'$ of retracts of $\CK$ in $D\CA T$. 
In particular, it is 
contained in the heart
$\CHQM$ of the motivic weight structure. The existence of weight 
filtrations~\ref{2A}~(4) for the weight structure $v$ then formally implies that
\[
D\CA T_{v \le 0} \subset \DgM_{F,w \le 0} \; ,
\]
and that
\[
D\CA T_{v \ge 0} \subset \DgM_{F,w \ge 0} \; .
\]
Now let $M_1 \in D\CA T_{w \le 0} = D\CA T \cap \DgM_{F,w \le 0}$.
Then for any $M_2 \in D\CA T_{v \ge 1}$, we have 
\[
\Hom_{D\CA T}(M_1,M_2) = 0 \; ,
\]
thanks to orthogonality~\ref{2A}~(3) for the motivic weight structure,
and to the fact that $D\CA T_{v \ge 1}$ is contained in 
$\DgM_{F,w \ge 1}$. Axioms~\ref{2A}~(1) and (4) easily imply
(see also \cite[Prop.~1.3.3~2]{Bo}) that $M_1 \in D\CA T_{v \le 0}$.
Therefore,
\[
D\CA T_{w \le 0} = D\CA T_{v \le 0} \; .
\]
In the same way, one proves that
\[
D\CA T_{w \ge 0} = D\CA T_{v \ge 0} \; .
\]
Altogether, the weight structure $v$ coincides with the data $w=w_\CA$. 
This proves part~(a) of our claim. We also see that 
part~(b) is formally implied by the following claim. 
(b')~The category $\CK$ is Abelian semi-simple.
(Since then $\CK$ will necessarily be pseudo-Abelian,
hence $D\CA T_{w = 0} = \CK'$ coincides with $\CK$.)   

Now consider two objects $N_1, N_2$ of $\CA$,
two integers $m_1,m_2$, and the group of morphisms
\[
\Hom \bigl( N_1(m_1)[2m_1] , N_2(m_2)[2m_2] \bigr) = 
\Hom \bigl( N_1 , N_2(m_2-m_1)[2(m_2-m_1)] \bigr)
\]
in $D \CA T$. Two essentially different cases occur: if $m_1 \ne m_2$,
then the group of morphisms is zero. Indeed, using descent
for finite extensions of $k$ as in the proof of Proposition~\ref{1H},
we reduce ourselves to the case $N_1 = N_2 = \BZ$,
where the desired vanishing follows from \cite[Prop.~4.2.9]{V}.

If $m_1 = m_2$, then
\[
\Hom (N_1(m_1)[2m_1] , N_2(m_2)[2m_2]) = \Hom (N_1 , N_2)
\]
can be calculated in the Abelian category $\CA$.

Thus in any of the two cases, the group
$\Hom (N_1(m_1)[2m_1] , N_2(m_2)[2m_2])$ coincides with 
\[
\Hom_{\Gr_\BZ \CA} \bigl( (N_1)_{m=m_1} , (N_2)_{m=m_2} \bigr) \; .
\]
Therefore, the functor defined in part~(c) of the claim
is fully faithful. Furthermore, it induces an equivalence
of categories between $\Gr_\BZ \CA$ and $\CK$.
The latter is therefore Abelian semi-simple.
This shows (b'), hence part~(b) of our claim. It also shows part~(c).
\end{Proof}

Statement (b) of Theorem~\ref{2E} should be considered as remarkable
in that it happens rarely that the heart of a given weight structure
is Abelian. We refer to \cite[Thm.~3.2]{P}, where this question is
studied abstractly. 

\begin{Cor} \label{2Ea}
The category $D \CA T$ is pseudo-Abelian. 
\end{Cor}

\begin{Proof}
By Theorem~\ref{2E}~(b), the heart
$D\CA T_{w = 0}$ is pseudo-Abelian and generates the triangulated
category $D\CA T$. Our claim thus follows from \cite[Lemma~5.2.1]{Bo}.
\end{Proof}

Here is another formal consequence of Theorem~\ref{2E}~(b).

\begin{Cor} \label{2Eb}
(a)~The inclusion of the heart 
$\iota_-: D \CA T_{w = 0} \into D \CA T_{w \le 0}$
admits a left adjoint 
\[
\Gr_0 : D \CA T_{w \le 0} \longto D \CA T_{w = 0} \; .
\]
For any $M \in D \CA T_{w \le 0}$, the adjunction morphism
$M \to \Gr_0 M$ gives rise to a weight filtration
\[
M_{\le -1} \longto M \longto \Gr_0 M \longto M_{\le -1}[1] 
\]
of $M$.
The composition $\Gr_0 \circ \iota_-$
equals the identity on $D \CA T_{w = 0}$. \\[0.1cm]
(b)~The inclusion of the heart 
$\iota_+: D \CA T_{w = 0} \into D \CA T_{w \ge 0}$
admits a right adjoint 
\[
\Gr_0 : D \CA T_{w \ge 0} \longto D \CA T_{w = 0} \; .
\]
For any $M \in D \CA T_{w \ge 0}$, the adjunction morphism
$\Gr_0 M \to M$ gives rise to a weight filtration
\[
\Gr_0 M \longto M \longto M_{\ge 1} \longto \Gr_0 M[1] 
\]
of $M$.
The composition $\Gr_0 \circ \iota_+$
equals the identity on $D \CA T_{w = 0}$.
\end{Cor}

\begin{Proof}
Let $M \in D \CA T_{w \le 0}$.
First choose an exact triangle
\[
M_{\le -2} \longto M \longto M_{-1,0} \longto M_{\le -2}[1] \; ,
\]
with $M_{\le -2} \in D \CA T_{w \le -2}$
and $M_{-1,0} \in D \CA T_{w \ge -1} \cap D \CA T_{w \le 0}$.
Orthogona\-lity~\ref{2A}~(3), together with the fact that
$M_{\le -2}[1] \in D \CA T_{w \le -1}$ 
shows that the morphism $M \to M_{-1,0}$
induces an isomorphism
\[
\Hom_{D \CA T} \bigl( M_{-1,0} , N \bigr) 
\isoto \Hom_{D \CA T} \bigl( M , N \bigr)
\]
for any object $N$ of the heart $D \CA T_{w = 0}$. Now choose an exact triangle
\[
M_{-1}' \longto M_{-1,0} \longto M_0' \stackrel{\alpha}{\longto} M_{-1}'[1] \; ,
\]
with $M_{-1}' \in D \CA T_{w = -1}$ (hence $M_{-1}'[1] \in D \CA T_{w = 0}$)
and $M_0' \in D \CA T_{w = 0}$. Recall that according to Theorem~\ref{2E}~(b),
$D \CA T_{w = 0}$ is Abelian semi-simple. Therefore,
the morphism $\alpha$ has a kernel and an image, 
both of which admit direct complements in $M_0'$ and in $M_{-1}'[1]$,
respectively. Choose a direct complement $M_0$ of $\ker \alpha$  in $M_0'$,
and a direct complement $M_{-1}[1]$ of $\imm \alpha$ in $M_{-1}'[1]$
(for some $M_{-1} \in D \CA T_{w = -1}$). 
\emph{Via} the restriction of $\alpha$, the object $M_0$ is isomorphic 
to the image. We thus get a commutative diagram
\[
\vcenter{\xymatrix@R-10pt{
        M_0' \ar@{>>}[d] \ar[r]^-{\alpha} &
        M_{-1}'[1] \ar@{>>}[d] \\
        \ker \alpha =: \Gr_0 M \ar[r]^-{0} &
        M_{-1}[1]
\\}}
\]
in $D \CA T_{w = 0}$ and in fact, a
morphism of exact triangles
\[
\vcenter{\xymatrix@R-10pt{
        M_{-1}' \ar[d] \ar[r] &
        M_{-1,0} \ar@{=}[d] \ar[r] &
        M_0' \ar[d] \ar[r]^-{\alpha} &
        M_{-1}'[1] \ar[d] \\
        M_{-1} \ar[r] &
        M_{-1,0} \ar[r] &
        \Gr_0 M \ar[r]^-{0} &
        M_{-1}[1]
\\}}
\]
in $D \CA T$.
By construction, and by orthogona\-lity~\ref{2A}~(3),
the morphism $M_{-1,0} \to \Gr_0 M$
induces an isomorphism
\[
\Hom_{D \CA T} \bigl( \Gr_0 M , N \bigr) 
\isoto \Hom_{D \CA T} \bigl( M_{-1,0} , N \bigr)
\]
for any object $N$ of the heart $D \CA T_{w = 0}$.
Choosing a cone of the composition $M \to M_{-1,0} \to \Gr_0 M$,
we get exact triangles
\[
M_{\le -1} \longto M \longto \Gr_0 M \longto M_{\le -1}[1] 
\]
and 
\[
M_{\le -2} \longto M_{\le -1} \longto M_{-1} \longto M_{\le -2}[1] \; .
\]
The second of the two exact triangles, together with
stability of $D \CA T_{w \le -1}$ under extensions
(cmp.~\cite[Prop.~1.3.3~3]{Bo}) shows that
$M_{\le -1}$ belongs to $D \CA T_{w \le -1}$.
Therefore, the first is a weight filtration of $M$. By construction,
the morphism $M \to \Gr_0 M$
induces an isomorphism
\[
\Hom_{D \CA T} \bigl( \Gr_0 M , N \bigr) 
\isoto \Hom_{D \CA T} \bigl( M , N \bigr)
\]
for any object $N$ of the heart $D \CA T_{w = 0}$.
From this property, it is easy to deduce the functorial behaviour of $\Gr_0 M$.

This proves part~(a) of the claim; the proof of part~(b) is dual.
\end{Proof}

\begin{Rem}
(a)~As the proof shows, Corollary~\ref{2Eb} remains true in the general context
of weight structures on triangulated categories,
whose heart is Abelian semi-simple. \\[0.1cm]
(b)~This more general version
of Corollary~\ref{2Eb} should be compared to \cite[Prop.~2.2]{W3}.
The conclusions on the existence of the adjoints $\Gr_0$ are the same.
On the one hand,
Corollary~\ref{2Eb} works without the additional assumption from \loccit \ on 
the absence of the adjacent weights $-1$ and $1$. On the other hand,
the fact that the heart is Abelian semi-simple, as we have just seen, is 
a vital ingredient of the proof.
There is another subtle difference between the two situations: In the setting of
\cite[Prop.~2.2~(a)]{W3}, the term $M_{\le -2}$ also behaves functorially
in $M$. This should not be expected to hold for
the term $M_{\le -1}$ from Corollary~\ref{2Eb}~(a).
\end{Rem}

\begin{Cor} \label{2Ec}
Let $M \in D \CA T_{w \ge -1} \cap D \CA T_{w \le 0}$.
Then the adjunction morphism
$M \to \Gr_0 M$ admits a right inverse and a kernel
(in the category $D \CA T$).
The latter is pure of weight $-1$.
Any choice of right inverse induces an isomorphism  
\[
M_{-1} \oplus \Gr_0 M \cong M 
\]
between $M$ and the direct sum of 
$M_{-1} \in D \CA T_{w = -1}$ and of $\Gr_0 M$.
\end{Cor}

\begin{Proof}
Either look at the proof of Corollary~\ref{2Eb}~(a). Or use its statement:
indeed, the adjunction morphism can be extended to a weight filtration
\[
M_{-1} \longto M \longto \Gr_0 M \stackrel{\alpha}{\longto} M_{-1}[1] 
\]
of $M$, with some $M_{-1} \in D \CA T_{w = -1}$. Since $\Gr_0$ is
left adjoint to $\iota_-$, and $M_{-1}[1] \in D \CA T_{w = 0}$,
the morphism $\alpha$ is necessarily zero. Therefore, the 
weight filtration splits.
\end{Proof}

Of course, the object $M_{-1}$ occurring in Corollary~\ref{2Ec}
is just the shift by $-1$ of $\Gr_0$ applied to $M[1] \in D \CA T_{w \ge 0}$.

\begin{Def}[{\cite[Def.~1.10]{W3}}] \label{2F}
Let $r \le s$ be two integers, and $D$ one of the categories
$D\CA T$ or $\DQgM$. An object $M$ of 
$D$ is said to be \emph{without weights $r,\ldots,s$} 
if there is an exact triangle
\[
M_{\le r-1} \longto M \longto M_{\ge s+1} \longto M_{\le r-1}[1]
\]
in $D$, with $M_{\le r-1} \in D_{w \le r-1}$ and
$M_{\ge s+1} \in D_{w \ge s+1}$.
\end{Def}
 
For the sequel, it will be important to know that the property 
of being without weights $r,\ldots,s$ is stable under extensions.
Recall that $L$ is said to be an extension of $M$ by $K$ if there
is an exact triangle
\[
K \longto L \longto M \longto K[1] \; .
\]

\begin{Prop} \label{2G}
Let 
\[
K \longto L \longto M \longto K[1] 
\]
be an exact triangle in $D\CA T$ or in $\DQgM$. 
Assume that $K$ and $M$ are both
without weights $r,\ldots,s$.
Then $L$ is 
without weights $r,\ldots,s$.
\end{Prop}

\begin{Proof}
According to the context, write $D$ for the category
$D\CA T$ resp.\ $\DQgM$ we are working in.
Let
\[
K_{\le r-1} \longto K \longto K_{\ge s+1} \longto K_{\le r-1}[1] 
\]
and
\[
M_{\le r-1} \longto M \longto M_{\ge s+1} \longto M_{\le r-1}[1] 
\]
be exact triangles in $D$, with
\[
K_{\le r-1} \, , \, M_{\le r-1} \in D_{w \le r-1}
\]
and
\[
K_{\ge s+1} \, \, , M_{\ge s+1} \in D_{w \ge s+1} \; . 
\]
By orthogonality~\ref{2A}~(3), there are no non-zero morphisms between
$M_{\le r-1}[-1]$ and $K_{\ge s+1}$. By \cite[Lemma~1.1]{L}
(with $f:Z_1 \to Z_2$ equal to the morphism $M[-1] \to K$), this implies
the existence of exact triangles
\[
M_{\le r-1}[-1] \longto K_{\le r-1} \longto L' \longto M_{\le r-1} \; ,
\]
\[
M_{\ge s+1}[-1] \longto K_{\ge s+1} \longto L'' \longto M_{\ge s+1} \; ,
\]
and 
\[
L' \longto L \longto L'' \longto L'[1] \; .
\]
The first of these triangles shows that $L' \in D_{w \le r-1}$.
The second shows that $L'' \in D_{w \ge s+1}$. 
Therefore, the third shows 
that $L$ is indeed without weights $r,\ldots,s$.
\end{Proof}

\begin{Rem}
As the proof shows, Proposition~\ref{2G} remains true in the general context
of weight structures on triangulated categories.
\end{Rem}


\bigskip


%
%

\section{The case of an algebraic base field: the $t$-structure}
\label{3}



In this section, we assume $k$ to be algebraic over the field $\BQ$
of rational numbers.
We first show that the data from Definition~\ref{1I}
define a $t$-structure on the triangulated category 
$D\CA T$ (Theorem~\ref{1Main}), and more generally,
on $D\CA T_{[a,b]}$ (Variant~\ref{1Var}). This provides a generalization
of the main result from \cite{L}
(which concerns the case of Tate motives).
Our strategy of proof is identical to the one from \loccit .
We then proceed 
(Theorem~\ref{2Main}) to give a
characterization of the weight structure on $D\CA T$
in terms of this $t$-structure. 
Specializing further to the case of number fields,
Theorem~\ref{2Main} implies a criterion allowing to identify
the weight structure \emph{via} the Hodge theoretic or
$\ell$-adic realization (Theorem~\ref{2H}). 

\begin{Thm} \label{1Main}
The pair $(D\CA T^{\le 0} , D\CA T^{\ge 0})$
(Definition~\ref{1I})
is a $t$-struc\-ture on $D\CA T$.
It has the following properties.
\begin{enumerate}
\item[(a)] The $t$-structure is non-degenerate. 
\item[(b)] Its heart $\CA T$ is generated (as a full Abelian sub-category
of $D\CA T$ stable under extensions) by the objects $N(m)$,
for $N \in \CA$ and $m \in \BZ$.
\item[(c)] Each object $M$ of $\CA T$ has a canonical 
\emph{weight filtration} by sub-objects
\[
0 \subset \ldots \subset W_{n-1} M \subset W_n M \subset 
\ldots \subset M \; .
\]   
This filtration is functorial and exact in $M$.
It is uniquely characterized by the properties of being finite
(\emph{i.e.}, $W_n M = 0$ for $n$ very small and $W_n M = M$ for $n$ very large), 
and of admitting sub-quotients 
\[
\gr_n \! M := W_n M / W_{n-1} M \; , \; n \in \BZ
\] 
of the form $N_n(-n/2)$, for some $N_n \in \CA$.
\item[(d)] The functor
\[
\bigoplus_{m \in \BZ} \gr_{2m}(m) :  \CA T \longto \Gr_\BZ \CA \; , \;
M \longmapsto \bigl( (\gr_{2m} \! M)(m) \bigr)_m
\]
is a faithful exact tensor functor to the $\BZ$-graded category over $\CA$. 
It thus identifies $\CA T$ with a tensor sub-category of $\Gr_\BZ \CA$.
\item[(e)] The natural maps
\[
\Ext^p_{\CA T} \bigl( M_1,M_2 \bigr) \longto 
\Hom_{D\CA T} \bigl( M_1,M_2[p] \bigr) 
\]
($\Ext^p = $ Yoneda Ext-group of $p$-extensions)
are isomorphisms, for all $p$, and all $M_1,M_2 \in \CA T$. Both sides
are zero for $p \ge 2$.
In particular, the Abelian category $\CA T$ is of cohomological
dimension one.
\end{enumerate}
\end{Thm}

We thus get in particular
the existence of two generating Abelian sub-categories, namely $\CA T$ and
$D\CA T_{w = 0}$, of the same triangulated category $D\CA T$.
The first of these is of cohomological dimension one, and the se\-cond
is semi-simple. In addition (Theorems~\ref{1Main}~(d)
and \ref{2E}~(c)), the first is
abstractly tensor equivalent to a tensor sub-category 
of the second. 

\begin{Rem}
In \cite[Def.~4.4.1]{Bo}, the notion of a $t$-structure \emph{adjacent}
to a given weight structure is defined. It may be important to note
that the $t$-structure $(D\CA T^{\le 0} , D\CA T^{\ge 0})$ 
from Theorem~\ref{1Main}
is \emph{not} adjacent to the motivic weight structure 
$(D\CA T_{w \le 0} , D\CA T_{w \ge 0})$ on $D\CA T$
studied in the previous section. Else, by definition we would have
\[
D\CA T_{w \le 0} = D\CA T^{\le 0} \quad \text{or} \quad
D\CA T_{w \ge 0} = D\CA T^{\ge 0}
\]
(according to whether we are in a situation of left or right adjointness).
But according to Theorem~\ref{2E}~(b), any object of the form
$\BZ(m)[2m]$, $m \in \BZ$, lies in the heart $D\CA T_{w = 0}$, while
$\BZ(-1)[-2] \in D\CA T^{\ge 2}$ and
$\BZ(1)[2] \in D\CA T^{\le -2}$.
\end{Rem}

Theorem~\ref{1Main} is the special case $(a,b) = (-\infty,\infty)$ of the
following.

\begin{Var}[cmp.\ {\cite[Thm.~1.4, Cor.~4.3]{L}}] \label{1Var}
Fix $a \le b$.
Then the pair
$(D\CA T_{[a,b]}^{\le 0} , D\CA T_{[a,b]}^{\ge 0})$ is a $t$-structure
on $D\CA T_{[a,b]}$.
It has the following properties.
\begin{enumerate}
\item[(a)] The $t$-structure is non-degenerate. 
\item[(b)] Its heart $\CA T_{[a,b]}$ is generated 
(as a full Abelian sub-category
of $D\CA T_{[a,b]}$ (or of $\DQgM$) stable under extensions) 
by the objects $N(-n/2)$,
for $N \in \CA$ and $a \le n \le b$.
\item[(c)] Each object $M$ of $\CA T_{[a,b]}$ has a canonical
weight filtration by sub-objects
\[
0 = W_{a-1} M \subset W_a M \subset 
\ldots \subset W_{b-1} M \subset W_b M = M \; .
\]   
This filtration is functorial and exact in $M$.
It is uniquely characterized by the property of 
admitting sub-quotients 
\[
W_n M / W_{n-1} M \; , \; n \in \BZ
\] 
of the form $N_n(-n/2)$, for some $N_n \in \CA$. For all $n \in \BZ$, we have
\[
W_n M / W_{n-1} M = \gr_n \! M 
\]  
as objects of the heart $\CA T_n$ of $D\CA T_n$.
\item[(d)] The functor
\[
\bigoplus_{m \in \BZ,a \le 2m \le b} \gr_{2m}(m) :  
\CA T_{[a,b]} \longto \bigoplus_{m \in \BZ,a \le 2m \le b} \CA 
\]
is a faithful exact tensor functor. 
\item[(e)] The natural maps
\[
\Ext^p_{\CA T_{[a,b]}} \bigl( M_1,M_2 \bigr) \longto 
\Hom \bigl( M_1,M_2[p] \bigr) 
\]
($\Hom =$ morphisms in $D\CA T_{[a,b]}$ (or in $\DQgM$))
are isomorphisms, for all $p$, and all $M_1,M_2 \in \CA T_{[a,b]}$. 
Both sides are zero for $p \ge 2$.
In particular, the Abelian category $\CA T_{[a,b]}$ is of cohomological
dimension one.
\item[(f)] For $a' \le a$ and $b \le b'$, the inclusion of 
$D\CA T_{[a,b]}$ into $D\CA T_{[a',b']}$ as a full triangulated
sub-category is compatible with the $t$-structures. That is, the
$t$-structure on $D\CA T_{[a,b]}$ is induced by the $t$-structure
on $D\CA T_{[a',b']}$. 
\end{enumerate}
\end{Var}

\begin{Proof}
The decisive ingredient is the following generalization of the vanishing
from \cite[Thm.~1.4]{L}:
\[
\Hom_{D\CA T_{[a,b]}} \bigl( N_1(m_1),N_2(m_2)[s] \bigr) = 0 \; , \;
\forall \, m_1 < m_2 \, , \, N_1, N_2 \in \CA \, , \, s \le 0 \; .
\]
It holds because $\Hom_{D\CA T_{[a,b]}} = \Hom_{\DQgM}$,
and $\Hom_{DM_{gm}(\bullet)_F}$ satisfies descent
for finite extensions $L / k$ of the base field. 
Choosing an extension $L$ splitting both 
$N_1$ and $N_2$ therefore allows to deduce the desired vanishing
from the Beilinson--Soul\'e vanishing conjecture 
\[
\Hom_{DM_{gm}(L)_\BQ} \bigl( \BZ(m_1),\BZ(m_2)[s] \bigr) \; ,
\]
which by the work of Borel is known for all number fields,
hence also for direct limits $L$ of such.

We now faithfully
imitate the proof of \cite[Thm.~1.4]{L}, to get assertions (a), (b), and (d).
We also get the following: the filtration $W_{\bullet} M$ induced by the
grading $\gr_\bullet M$ is functorial. By construction, the sub-quotient
$\gr_n M$ lies in $\CA T_n$.
Its unicity follows from the fact that there are no non-zero
morphisms from objects of $\CA T$ of weights at most $r$ to objects
of weights at least $r+1$. To prove this, use induction on the lenght
of weight filtrations, and the vanishing 
\[
\Hom_{D\CA T} \bigl( N_1(m_1)[r],N_2(m_2)[s] \bigr) = 0 \; , \;
\forall \, m_1 > m_2 \, , \, N_1, N_2 \in \CA \, , \, r,s \in \BZ 
\]
(see the proof of Proposition~\ref{1H}). 
We thus get part~(c) of our claim.

Part (f) follows from the definition of our $t$-structure, and from the
compatibility of the functors $\gr_n$ under the inclusion of
$D\CA T_{[a,b]}$ into $D\CA T_{[a',b']}$.

As for claim~(e), we faithfully imitate the proof of 
\cite[Cor.~4.3]{L}.
\end{Proof}

\begin{Cor} \label{1N}
The identity on $\CA T$ extends canonically to
an equivalence of triangulated categories 
\[
D^b \bigl( \CA T) \longto D \CA T 
\]
between the bounded derived category of $\CA T$ and $D \CA T$.
Its composition with the cohomology functor
$D \CA T \to \CA T$ associated to the $t$-structure 
of Theorem~\ref{1Main} equals the canonical cohomology
functor on $D^b ( \CA T )$. 
\end{Cor}

\begin{Proof}
Recall the definition of the category $\ShN$ 
of \emph{Nisnevich sheaves with transfers}
\cite[Def.~3.1.1]{V}.  
It is Abelian \cite[Thm.~3.1.4]{V}, and there is a canonical full
triangulated embedding 
\[
\DeffgM \longinto D^- \bigl( \ShN \bigr)
\]
into the derived category of complexes of Nisnevich sheaves
bounded from above \cite[Thm.~3.2.6, p.~205]{V}.
Imitating the construction from \loccit \ using $F$ instead
of $\BZ$ as ring of coefficients,
one shows that there is a canonical full
triangulated embedding 
\[
\DeffQgM \longinto D^- \bigl( \ShN_F \bigr) \; ,
\]
where $\ShN_F$ denotes the  Abelian cate\-go\-ry
of Nisnevich sheaves with transfers taking values in $F$-modules.
We thus get a canonical embedding
into $D \bigl( \ShN_F \bigr)$ of any full triangulated category
$\CC$ of $\DeffQgM$, and hence in particular
for $\CC = D \CA T$. Our claim thus follows from 
\cite[Thm.~1.1~(a), (d)]{W4}:
indeed, 
$\Hom_{D \CA T} (M_1,M_2[2]) = 0$ for any two objects
$M_1, M_2$ in $\CA T$ (Theorem~\ref{1Main}~(e)), and $\CA T$ generates $D \CA T$
(Theorem~\ref{1Main}~(b)). 
\end{Proof}

We already mentioned the special cases $\CA = triv$ and
$\CA = \QAM$. A third case appears worthwhile mentioning.

\begin{Def} \label{1K}
(a)~Define the category $\QDM$ as the full Abelian $F$-linear sub-category of
$\QAM$ of objects on which the Galois group acts \emph{via} a commutative (finite)
quotient. \\[0.1cm]
(b)~Define the \emph{triangulated
category of Dirichlet--Tate motives over $k$}
as the strict full tensor triangulated sub-category $\DQDTM$ of $\DQgM$ 
generated by $\QDM$ and $\DQTM$. 
\end{Def}

Similarly, for any algebraic extension
$K$ of $k$, we could define the \emph{triangulated category of
Artin-Tate} (resp.\ \emph{Dirichlet--Tate}, resp...)
\emph{motives over $k$ trivializable over $K$} 
by letting $\CA$ equal the full Abelian $F$-linear 
sub-category $MA(K/k)_F$ (resp.\ $MD(K/k)_F$, resp...)
of $\QAM$ (resp.\ $\QDM$, resp...) of objects on which the absolute
Galois group of $K$, when identified with a subgroup of the
Galois group of $k$, acts trivially. 

\begin{Cor} \label{1L}
The conclusions of Theorem~\ref{1Main}, Variant~\ref{1Var}
and Corollary~\ref{1N} hold in particular in any of the following three cases.
\begin{enumerate}
\item[(1)] $\CA = triv$. In particular, 
this gives back the main result of \cite{L}.
The heart $\CA T$ equals the Abelian category 
$MT(k)_F$ of \emph{mixed Tate motives}.
\item[(2)] $\CA = \QDM$. In this case, the category $D \CA T$
equals the triangulated category $\DQDTM$ of Dirichlet--Tate motives.
Its heart $\CA T$ equals the Abelian category $\QDTM$ of
\emph{mixed Dirichlet--Tate motives}.
\item[(3)] $\CA = \QAM$. In this case, the category $D \CA T$
equals the triangulated category $\DQATM$ of Artin--Tate motives.
Its heart $\CA T$ equals the Abelian category $\QATM$ of
\emph{mixed Artin--Tate motives}. 
\end{enumerate}
\end{Cor}

\begin{Rem}
(a)~An equivalent construction of the 
category $\QATM$, for $F = \BQ \,$, is given in \cite[Sect.~2.17]{DG}. \\[0.1cm]
(b)~Note that by construction, an inclusion $\CA \subset \CB$
of strict
full Abelian semi-simple $F$-linear tensor sub-categories of $\QAM$
containing $triv$
induces first a strict full tensor triangulated embedding 
$D \CA T \subset D \CB T$,
and then a strict full exact tensor embedding $\CA T \subset \CB T$.
An object of $D \CB T$ belongs to $D \CA T$ if and only if its cohomology
objects (with respect to the $t$-structure from Theorem~\ref{1Main}) 
lie in $\CA T$. The equivalences of Corollary~\ref{1N} for $\CA$ and $\CB$
fit into a commutative diagram\[
\vcenter{\xymatrix@R-10pt{
        D^b (\CA T) \ar[d] \ar[r]^-{\cong} &
        D \CA T \ar@{^{ (}->}[d] \\
        D^b (\CB T) \ar[r]^-{\cong} &
        D \CB T
\\}}
\]
\cite[Thm.~1.1~(b)]{W4}.
In particular, the bounded derived category $D^b (\CA T)$ 
is canonically identified with a full sub-category of 
the bounded derived category $D^b (\CB T)$.
\end{Rem}

\begin{Def} \label{1M}
Let $M$ be a mixed Artin--Tate motive, with weight filtration
\[
0 \subset \ldots \subset W_{r-1} M \subset W_r M \subset 
\ldots \subset M \; .
\] 
Let $n$ be an integer. \\[0.1cm]
(a)~We say that $M$ is \emph{of weights $\le n$} if $W_n M = M$. \\[0.1cm]
(b)~We say that $M$ is \emph{of weights $\ge n$} if $W_{n-1} M = 0$. \\[0.1cm]
(c)~We say that $M$ is \emph{pure of weight $n$} if it is both of
weights $\le n$ and of weights $\ge n$, \emph{i.e.}, if $W_{n-1} M = 0$
and $W_n = M$. \\[0.1cm]
(d)~We say that $M$ is \emph{without weight $n$}
if $W_{n-1} M = W_n M$, \emph{i.e.}, if the sub-quotient
$W_n M / W_{n-1} M$ is trivial. 
\end{Def}

Of course, any mixed Artin--Tate motive is without weight $n$,
whenever $n$ is odd.
Denote by 
\[
\tau^{\le n} \; , \; \tau^{\ge n} : D\CA T \longto D\CA T
\]
the truncation functors, and by 
\[
\CH^n : D\CA T \longto \CA T 
\]
the cohomology functors associated to the $t$-structure from
Theorem~\ref{1Main}. 

\begin{Thm} \label{2Main} 
Let $K \in D\CA T$, and $r \le s$. \\[0.1cm]
(a)~$K$ lies in the heart $D\CA T_{w=0}$ of $w$ if and only if 
the object $\CH^n \! K$ of $\CA T$ is pure of weight $n$,
for all $n \in \BZ$. \\[0.1cm]
(b)~$K$ lies in $D\CA T_{w \le r}$ if and only if 
$\CH^n \! K$ is of weights $\le n + r$, 
for all $n \in \BZ$. \\[0.1cm] 
(c)~$K$ lies in $D\CA T_{w \ge s}$ if and only if 
$\CH^n \! K$ is of weights $\ge n + s$, 
for all $n \in \BZ$. \\[0.1cm]
(d)~$K$ is without weights $r,\ldots,s$ if and only if 
$\CH^n \! K$ is without weights $n + r,\ldots,n + s$, 
for all $n \in \BZ$. 
\end{Thm}

\begin{Proof}
Observe that the triangulated category $D\CA T$ is generated 
by the heart $D\CA T_{w=0}$ of $w$ (Theorem~\ref{2E}~(b))
as well as by the heart $\CA T$ of $t$ (Theorem~\ref{1Main}~(a)).
This will allow to simplify the proof.

The explicit description of objects $K$ of $D\CA T_{w=0}$
from Theorem~\ref{2E}~(b) shows that the
$\CH^n \! K$ are indeed pure of weight $n$, for all $n$
(see Theorem~\ref{1Main}~(c)). To show that any $K$ whose cohomology
objects $\CH^n \! K$ are pure of weight $n$, does belong to $D\CA T_{w=0}$,
we may assume by the above 
that $K$ is concentrated in one degree (with respect
to the $t$-structure), say $K = M[d]$ for some $M \in \CA T$ and $d \in \BZ$.
By assumption, the mixed Artin-Tate motive $M$ is pure of weight $-d$,
and hence (Theorem~\ref{1Main}~(c)) of the form $N(d/2)$, for some Artin
motive $N$ belonging to $\CA$. The latter is clearly a Chow motive, and
hence so is its tensor product with the Chow motive $\BZ(d/2)[d]$.
Therefore, $K$ is a Chow motive belonging to $D\CA T$.
By Theorem~\ref{2E}~(b), it is in the heart $D\CA T_{w=0}$. 
This shows part~(a). 

We leave it to the reader to deduce (b) and (c) from (a).  

As for part~(d), it is easy to see that the cohomology 
$\CH^n \! K$ of an object $K \in D\CA T$ 
without weights $r,\ldots,s$ is without weights $n+r,\ldots,n+s$,
for all $n$
(use (b) resp.\ (c) for the constituents of a suitable weight filtration
of $K$). To prove the inverse implication, we use induction on the number
of integers $n$ such that $\CH^n \! K \ne 0$. 
If this number equals one, then $K = M[d]$ for some $M \in \CA T$ and 
$d \in \BZ$.
By assumption, the mixed Artin-Tate motive $M$ is without weights 
$-d+r,\ldots,-d+s$. By Definition~\ref{1M}~(d), its weight filtration
thus satisfies the relation
\[
W_{-d+r-1} M = W_{-d+r} M = \ldots = W_{-d+s} M \; .
\]
The sequence
\[
0 \longto W_{-d+r-1} M \longto M \longto M / W_{-d+s} M \longto 0
\]
is therefore exact in $\CA T$. It gives rise to an exact triangle
\[
\bigl( W_{-d+r-1} M \bigr)[d] \longto K \longto 
\bigl( M / W_{-d+s} M \bigr)[d] \longto \bigl( W_{-d+r-1} M \bigr)[d+1]
\]
in $D\CA T$. By parts (b) and (c),
\[
\bigl( W_{-d+r-1} M \bigr)[d] \in D\CA T_{w \le r-1} \; ,
\]
and
\[
\bigl( M / W_{-d+s} M \bigr)[d] \in D\CA T_{w \ge s+1} \; .
\]
Therefore, the object $K$ is indeed without weights $r,\ldots,s$.

By Proposition~\ref{2G}, the property of being without weights $r,\ldots,s$
is stable under extensions in $D\CA T$. 
This allows to perform the induction step.
\end{Proof}

To conclude, let us now consider realizations
(\cite[Sect.~2.3 and Corrigendum]{H}; see \cite[Sect.~1.5]{DG}
for a simplification of this approach). We  
assume from now on that $k$ is a number field, and concentrate on two
realizations (the statement from Theorem~\ref{2H}
below then formally generalizes to
any of the other realizations ``with weights'' considered in \cite{H}):
\begin{enumerate}
\item[(i)] the Hodge theoretic realization
\[
R_\sigma : \DQgM \longto D
\]
associated to a fixed embedding $\sigma$ of the number field $k$ into
the field $\BC$ of complex numbers. Here, $D$ is the bounded derived category 
of mixed graded-polarizable $\BQ$-Hodge structures 
\cite[Def.~3.9, Lemma~3.11]{Be}, tensored with $F$,
\item[(ii)] the $\ell$-adic realization
\[
R_{\ell} : \DQgM \longto D
\]
for a prime $\ell$. Here, $D$ is the bounded ``derived category'' of
constructible $\BQ_\ell$-sheaves on $\Spec (k)$
\cite[Sect.~6]{E}, tensored with $F$.
\end{enumerate}

Choose and fix one of these two, denote it by $R$,
recall that it is a contravariant tensor functor,
and use the same letter for its restriction to the
sub-category $D\CA T$ of $\DQgM$. 
The category $D\CA T$ 
is equipped with a $t$-structure. The same is true for $D$;
write $H^n$ for the cohomology functors.
It is easy to see that $R$ is $t$-exact (since it maps $\CA T$ to the
heart of $D$).
In particular, it induces an exact contravariant functor $R_0$ from
the heart $\CA T$ of $D\CA T$ to the heart of $D$, which we shall
denote by $\CB$. 
As for the weight structure on $D\CA T$, note that $R_0$ maps the
pure Tate motive $\BZ(m)$ to the pure Hodge structure $\BQ(-m)$
(when $R=R_\sigma$) and to the pure $\BQ_\ell$-sheaf $\BQ_\ell(-m)$
(when $R=R_\ell$), respectively \cite[Thm.~2.3.3]{H}. The following is a consequence of 
Theorem~\ref{1Main}.

\begin{Cor} \label{2Ha}
Assume $k$ to be a number field. \\[0.1cm]
(a)~The realization
\[
R : D\CA T \longto D
\]
is conservative. In other words, an object $K$ of $D\CA T$ is zero
if and only if its image $R(K)$ under $R$ is. \\[0.1cm]
(b)~The induced functor
\[
R_0 : \CA T \longto \CB
\]
is conservative. \\[0.1cm]
(c)~The functor $R_0$ respects and detects
weights up to inversion of the sign.
More precisely,
an object $M$ of $\CA T$ is pure of weight $n$ if and only if $R_0(M)$
is pure of weight $-n$.
\end{Cor}

Note that there is a notion of purity and mixedness for objects of $\CB$. \\

\begin{Proofof}{Corollary~\ref{2Ha}}
Let $K \in D\CA T$. Given the $t$-exactness and contravariance of $R$,
we have the formula 
\[
H^n \! R(K) = R_0 \bigl( \CH^{-n} \! K \bigr) 
\]
for all $n$. By Theorem~\ref{1Main}~(a),
the $t$-structure on $D\CA T$ is non-degenerate.
Hence (b) implies (a).

Recall that by Theorem~\ref{1Main}~(c),
there is a unique finite weight filtration on any object of $\CA T$.
Also, the functor $R_0$ is exact.
Hence (b) is implied by conservativity of the restriction of $R_0$ 
to the sub-category of objects of $\CA T$, which are pure of some weight.
But this property is
clearly implied by (c) (since the zero object of $\CB$ 
is pure of any weight).

Let $M \in \CA T$. As before,
we may assume that
$M$ is pure of some weight, say $n$. Again by 
Theorem~\ref{1Main}~(c), $M$ is of the form $N(-n/2)$, 
for some Artin-Tate motive $N$. Thus, $R_0(M) \cong R_0(N)(n/2)$
is pure of weight $-n$. It is zero if and only if $R_0(N)$ is,
which is the case if and only if $N$ is. 
\end{Proofof}

We finally get the characterization of the weight structure we aimed at.

\begin{Thm} \label{2H}
Assume $k$ to be a number field. 
Then the realization $R$ respects and
detects the weight structure. More precisely,
let $K \in D\CA T$, and $r \le s$. \\[0.1cm]
(a)~$K$ lies in the heart $D\CA T_{w=0}$ of $w$ if and only if 
the $n$-th cohomology object $H^n \! R(K) \in \CB $ of $R(K)$ 
is pure of weight $n$, 
for all $n \in \BZ$. \\[0.1cm]
(b)~$K$ lies in $D\CA T_{w \le r}$ if and only if 
$H^n \! R(K)$ is of weights $\ge n - r$, 
for all $n \in \BZ$. \\[0.1cm] 
(c)~$K$ lies in $D\CA T_{w \ge s}$ if and only if 
$H^n \! R(K)$ is of weights $\le n - s$, 
for all $n \in \BZ$. \\[0.1cm]
(d)~$K$ is without weights $r,\ldots,s$ if and only if 
$H^n \! R(K)$ is without weights $n - s,\ldots,n - r$, 
for all $n \in \BZ$.  
\end{Thm}

\begin{Proof}
Recall that
$H^n \! R(K) = R_0 \bigl( \CH^{-n} \! K \bigr)$.
The claim thus follows from 
Theorem~\ref{2Main} and Corollary~\ref{2Ha}.
\end{Proof}

\begin{Rem}
As the proof shows, the analogues of parts (a) and (b) of
Corollary~\ref{2Ha} continue to hold for any of the realizations
(including those ``without weights'') considered in \cite{H}. 
This is true in particular for 
\begin{enumerate}
\item[(iii)] the Betti realization, \emph{i.e.}, the composition of the
Hodge theoretic realization $R_\sigma$ with the forgetful functor
to the bounded derived category of $F$-modules of finite type,
\item[(iv)] the topological $\ell$-adic realization, \emph{i.e.}, the
composition of the $\ell$-adic realization $R_{\ell}$
with the forgetful functor to the bounded derived category of
$F \otimes_\BQ \BQ_\ell$-modules of finite type \cite[Thm.~7.2~(i)]{E}.
\end{enumerate}
\end{Rem}

\begin{Rem} \label{2I}
(a)~For the Hodge theoretic realization
\[
R = R_\sigma : \DQgM \longto D 
\]
($D= $ the bounded derived category 
of mixed graded-polarizable $\BQ$-Hodge structures, 
tensored with $F$),
it is possible to give a more conceptual interpretation of 
respect of the weight structure. In fact, there is a canonical weight structure
$w_\FH$ on $D$, characterized by the property of admitting as heart 
the full sub-category $\CK$ of classes
of complexes $K$ of Hodge structures, whose $n$-th cohomology object 
$H^n \! K$ is pure of weight $n$, 
for all $n \in \BZ$. In order to prove this claim,
apply \cite[Thm.~4.3.2~II~1 and 2]{Bo}, observing that
(1)~$\CK$ generates the triangulated category $D$,
(2)~$\CK$ is negative: 
\[
\Hom_D(M_1,M_2[i]) = 0
\]
for any two objects $M_1, M_2$ of $\CK$, and any integer $i > 0$
(it is here that the polarizability assumption enters), and (3)~any retract of 
an object of $\CK$ in $D$ belongs already to $\CK$. \\

To say that $R : \DQgM \to D$ 
respects the weight structure means then that $R$ respects the 
pairs of sub-categories $(\DgM_{F,w \le 0} , \DgM_{F,w \ge 0})$
and $(D_{w \le 0} , D_{w \ge 0})$:
\[
R \bigl( \DgM_{F,w \le 0} \bigr) \subset D_{w \ge 0} \quad , \quad
R \bigl( \DgM_{F,w \ge 0} \bigr) \subset D_{w \le 0} 
\]
(recall that $R$ is contravariant). Given that 
$\DQgM$ is generated by its heart, this requirement is equivalent
to saying that $R$ respects the hearts,
i.e., that it maps $\CHQM$ to $\CK$ --- which is a true statement,
since the Hodge structure on the $n$-th Betti cohomology of a 
proper smooth variety is indeed pure of weight $n$, 
for all $n \in \BZ$. This observation implies immediately
the ``only if'' part of the statements
of Theorem~\ref{2H}. \\[0.1cm]
(b)~By contrast, for the $\ell$-adic realization
\[
R = R_{\ell} : \DQgM \longto D
\]
($D= $ the bounded ``derived category'' of
constructible $\BQ_\ell$-sheaves on $\Spec k$, tensored with $F$),
there is no such interpretation, since there is no reasonable 
weight structure on $D$. Indeed, according to \cite[Rem.~6.8.4~i)]{J},
for any odd integer $m \in \BZ$,
\[
\Hom_D(\BQ_\ell(0),\BQ_\ell(m)[1]) \ne 0 \; . 
\]
This is true in particular when $m$ is negative, i.e.,
$\BQ_\ell(m)[1]$ is pure of strictly positive weight $-2m+1$.
Therefore, orthogonality~\ref{2A}~(3) is violated. 
\end{Rem}


\bigskip

%
%

\end{document}